\documentclass[11pt]{amsart}
\usepackage[T1]{fontenc}
\usepackage{lmodern}
\usepackage[margin=1.1in]{geometry}
\usepackage{amsmath,amssymb,amsthm,mathtools}
\usepackage{microtype}
\usepackage{xcolor}
\usepackage{aliascnt}
\usepackage[colorlinks=true,linkcolor=blue!50!black,citecolor=blue!50!black,urlcolor=blue!50!black]{hyperref}
\usepackage[capitalize,nameinlink]{cleveref}

\newcommand{\C}{\mathbb C}

\newcommand{\Z}{\mathbb Z}

\newcommand{\II}{\mathrm{II}_1}
\newcommand{\E}{E_N}                          
\newcommand{\Tr}{\operatorname{Tr}}           
\newcommand{\Proj}{\operatorname{Proj}}
\newcommand{\U}{\mathcal U}
\newcommand{\Hh}{\mathcal H}
\newcommand{\sa}{\mathrm{sa}}

\newcommand{\one}{\mathbf 1}                  
\newcommand{\Lam}{\Lambda}                    
\newcommand{\basic}[2]{\langle #1,e_{#2}\rangle}

\newcommand{\norm}[1]{\left\lVert#1\right\rVert}
\newcommand{\abs}[1]{\left\lvert#1\right\rvert}

\newif\ifshowtodos

\showtodosfalse  

\newcommand{\newsharedtheorem}[3][theorem]{%
  \newaliascnt{#2}{#1}\newtheorem{#2}[#2]{#3}\aliascntresetthe{#2}}
\newtheorem{theorem}{Theorem}[section]
\newsharedtheorem{lemma}{Lemma}
\newsharedtheorem{proposition}{Proposition}
\newsharedtheorem{corollary}{Corollary}
\newtheorem{mainthm}{Theorem}

\newsharedtheorem[mainthm]{maincor}{Corollary}
\theoremstyle{definition}
\newsharedtheorem{definition}{Definition}
\newsharedtheorem{example}{Example}
\theoremstyle{remark}
\newsharedtheorem{remark}{Remark}
\numberwithin{equation}{section}

\crefname{mainthm}{Theorem}{Theorems}
\crefname{maincor}{Corollary}{Corollaries}

\title[The Carpenter property for subfactors]%
{The Carpenter property for subfactors\\of uniform infinite index}
\author{Dietmar Bisch}
\email{dietmar.bisch@vanderbilt.edu}
\author{Julio Cáceres}
\email{julio.e.caceres.gonzales@vanderbilt.edu}
\date{October 30, 2026}

\begin{document}

\begin{abstract}
We show that if $N\subset M$ is an inclusion of $\II$ factors of uniform
infinite index, for instance an irreducible inclusion of infinite
index, then every positive contraction $b\in N$ is of the form $E_N(p)$
for a projection $p\in M$. This gives the exact version of Popa's
asymptotic result in the ultrapower, and it is a subfactor analogue of
the Carpenter problem for masas. In particular, Popa's relative
dimension set of such an inclusion is $[0,1]$, so the set of
expectations of projections does not distinguish inclusions without
corners of finite index. As applications, the Carpenter property holds
for $N\subset N\rtimes G$ whenever an infinite group $G$ acts outerly
on a $\II$ factor $N$, and for ICC groups
$H\le G$ the inclusion $L(H)\subset L(G)$ has it exactly when
$[G:H]=\infty$, confirming a recent conjecture. We also give examples
showing that the Carpenter property does not imply uniform infinite
index, and that infinite index does not imply the Carpenter property.
\end{abstract}

\maketitle

\section{Introduction}\label{sec:intro}

A basic problem in the theory of subfactors is to find invariants of an
inclusion $N\subset M$ of $\II$ factors. In finite index the central
invariants are Jones' index \cite{Jones83} and the standard invariant,
the lattice of higher relative commutants, which is a complete
invariant for amenable inclusions of the hyperfinite factor
\cite{Popa94}. In infinite
index both lose most of their force: the index is $\infty$, the basic
construction is no longer a $\II$ factor, and the higher relative
commutants need not be finite-dimensional. Invariants of a different
kind are needed. Popa proposed to use the ways in which subalgebras of
$M$ can sit orthogonally to $N$ \cite[\S3.2]{Popa19}, in the sense of
\cite{Popa83}. The simplest of these are the expectations of
projections: the set $\E(\Proj(M))$, which Popa studies as a dilation
problem \cite[\S7]{Popa89}, \cite[\S3.6]{Popa19}, and its scalar part, the relative
dimension set $\Lam(M,N)$ of values $\alpha$ with $\E(p)=\alpha1$
\cite{Popa89}. We expected the \emph{Carpenter property}, namely that
$\E(\Proj(M))$ is the whole positive part of the unit ball of $N$, to
distinguish inclusions of infinite index that the classical invariants
cannot tell apart.

In finite index the Carpenter property always fails
\cite[Prop.~3.14]{HKS}. By the
Pimsner--Popa inequality $\E(x)\ge[M:N]^{-1}x$ for $x\in M_+$
\cite{PP86}, a nonzero projection $p$ with $\E(p)=\alpha1$ satisfies
$\alpha p=p\E(p)p\ge[M:N]^{-1}p$, so no scalar $0<\alpha<[M:N]^{-1}$
is the expectation of a projection (\cref{lem:PP}). In infinite index
this obstruction disappears. For inclusions of \emph{uniform infinite
index}, meaning that no corner $Np\subset pMp$ with $p\in N'\cap M$ has
finite index (equivalently, $M\nprec_MN$ in the sense of Popa's
intertwining by bimodules \cite{Popa06Rigidity}; see \cref{lem:WM}),
Popa proved that every positive contraction of $N$ is the expectation
of a projection after passing to the ultrapower: in
$N^\omega\subset M^\omega$, or equivalently up to arbitrarily small
error in $M$ \cite[\S3.6]{Popa19}. For irreducible inclusions, the
scalar case was already observed in \cite{Popa89}; see
\cite[\S3.3]{Popa19}. This leaves the exact question in $M$ itself.

The problem of studying the Carpenter property for general inclusions
$N\subset M$, where $M$ is a $\II$ factor, was initiated by Popa in
\cite[\S7]{Popa89}, who singled out the cases where $N$ is a maximal
abelian subalgebra (masa) or a subfactor of $M$. It was called the
Carpenter problem when $N$ is a masa in
\cite{Kadison02a,Kadison02b}, a terminology that is suggestive and has
already been adopted in the literature. This is the case of projections
in the Schur--Horn problem of \cite{AK06}.
Partial results were obtained in
\cite{ArgeramiMassey07,ArgeramiMassey09,BhatRavichandran,DFHS12}, and a
proof for every masa in a $\II$ factor was given in
\cite[Thm.~5.7]{Ravichandran}. The problem for general subalgebras was
also considered in \cite{HKS}, from the point of view of Lyapunov's convexity
theorem \cite{AkemannAnderson91}, where the case of group algebras was
conjectured. These methods use the abelian structure of a masa, a large
relative commutant $N'\cap M$, or the non-injectivity of $\E$ on every
corner of $M$ (\cref{rem:extreme}). Under uniform infinite index the
first two are not available in general, and we do not know whether the
third holds, so they do not apply here.

Our interest in the question began with crossed products
$N\subset N\rtimes G$, where the canonical unitaries form an infinite
orthonormal basis of $L^2(N\rtimes G)$ over $N$. The main result shows
that the Carpenter property only sees corners of finite index: once no
corner has finite index, every positive contraction of $N$ is the
expectation of a projection.

\subsection*{Main results}
\begin{mainthm}\label{thm:A}
Let $N\subset M$ be an inclusion of $\II$ factors of uniform infinite
index (\cref{def:uniform}); for instance, $N'\cap M=\C1$ and
$[M:N]=\infty$. Then $N\subset M$ has the \emph{Carpenter property}
(\cref{def:carpenter}): for every $b\in N$ with $0\le b\le1$ there is a
projection $p\in M$ with $\E(p)=b$. In particular, Popa's relative
dimension set is $\Lam(M,N)=[0,1]$ (\cref{def:Lambda}).
\end{mainthm}

The proof relies on a flattening theorem for completely positive maps
\cite[Cor.~1.3]{Jeong26} (\cref{thm:flattening}). No separability is
assumed; by contrast, in \cite[Prop.~3.10]{HKS} it is the
nonseparability of $M$ that produces the Carpenter property.

In a companion paper, the second-named author shows that $\Lam(M,N)$ is
also determined by the index for every extremal inclusion of finite
index \cite{Crelative}. Together with \cref{thm:A}, this shows that
$\Lam(M,N)$ depends only on
$[M:N]$ for every irreducible inclusion of $\II$ factors.

Uniform infinite index is easy to check when $L^2(M)$ has an infinite
orthonormal basis over $N$ made of unitaries of $M$
(\cref{prop:unitary-basis}). This covers crossed products and group
algebras, and proves the conjecture of \cite{HKS}.

\begin{maincor}\label{cor:B}
\begin{enumerate}
\item If $N$ is a $\II$ factor and $G$ is an infinite discrete group
 acting on $N$ such that $N\rtimes G$ is a factor, then
 $N\subset N\rtimes G$ has the Carpenter property.
\item If $H\le G$ are ICC groups, then $L(H)\subset L(G)$ has the
 Carpenter property if and only if $[G:H]=\infty$
 \cite[Conj.~3.15]{HKS}.
\end{enumerate}
\end{maincor}

Uniform infinite index cannot be replaced by infinite index.
\Cref{sec:examples} gives an inclusion of infinite index without the
Carpenter property, whose relative dimension set misses an interval,
and an inclusion with the Carpenter property that does not have uniform
infinite index. Thus the implications
\begin{equation}\label{eq:implications}
 \text{uniform infinite index}\ \Longrightarrow\
 \text{Carpenter property}\ \Longrightarrow\ [M:N]=\infty
\end{equation}
are both strict. For inclusions of infinite index, the Carpenter
property can therefore only detect corners of finite index, and the
examples show that such corners are detected in some cases and not in
others.

\subsection*{Strategy of the proof}
Popa's ultrapower result gives, for each $\epsilon>0$, a projection
$p_\epsilon\in M$ with $\norm{\E(p_\epsilon)-b}_2<\epsilon$, but nothing
controls the distance between these approximations, so they need not
converge. Another route is to show that the extreme points of the
convex set $\{x\in M:0\le x\le1,\ \E(x)=b\}$ are projections; this
holds when $\E$ is not injective on any corner of $M$
(\cref{rem:extreme}), which we do not know under uniform infinite index.
Instead, we improve an approximate lift while controlling how far it
moves.

For $b\in N$ with $\eta1\le b\le(1-\eta)1$, we
construct projections $p_k$ with $\E(p_k)\to b$ fast enough that
$(p_k)$ converges to a projection $p$ with $\E(p)=b$. Each step has two
parts. A convex duality argument, based on Sion's minimax theorem,
gives a positive contraction $x$ with $\E(x)=b$ whose affine cost
$\ell_p(x)=\tau(p)+\tau(x)-2\tau(px)$ is at most
$\eta^{-1/2}\norm{b-\E(p)}_1$ (\cref{lem:correction}). On projections
$q$ the cost is $\ell_p(q)=\norm{q-p}_2^2$, so it controls the
distance moved, provided $x$ can be replaced by a projection $q$ with
nearly the same expectation and nearly the same overlap $\tau(pq)$.
This is a flattening problem for the completely positive map
$y\mapsto(\E(y),\tau(py))$ on the spectral corners $eMe$ of $x$.
The flattening theorem of \cite{Jeong26} solves it once the bimodule
${}_{eMe}eL^2(M)_N$ is weakly mixing, and this is where uniform
infinite index enters (\cref{lem:uniformcorners}). The projection $q$
need not be close to $x$; only its expectation and its overlap with $p$
matter (\cref{prop:controlled}). Finally, a spectral decomposition of a
general $b$, together with the stability of uniform infinite index
under compression by projections of $N$ (\cref{lem:uniformNcorner}),
removes the spectral gap assumption.

\subsection*{Organization}
\Cref{sec:prelim} collects preliminaries: the basic construction, the
relative dimension set, the Carpenter property and uniform infinite
index. \Cref{sec:correction} gives the convex correction,
\cref{sec:flattening} proves weak mixing and flattening on corners, and
\cref{sec:lifting} approximates by projections and constructs the exact
lift, proving \cref{thm:A}. \Cref{sec:bases} treats inclusions with a
unitary orthonormal basis and proves \cref{cor:B}.
\Cref{sec:examples} gives examples without uniform infinite index, and
\cref{sec:final}
collects further remarks and open problems.

\subsection*{AI use statement}
Most of the proofs and arguments in this paper were developed with the
help of ChatGPT-6 Astra (OpenAI), which worked under the authors'
guidance. The
writing and organization of the paper, including the checking of
references against the cited sources, were carried out with Claude
Opus 5.5 (Anthropic), with editorial direction from the authors. The
authors directed the project, made the mathematical and editorial
decisions, checked all proofs, and take full responsibility for the
content of the paper.

\subsection*{Acknowledgements}
D.B. and J.C. were partially supported by U.S. Army Research Office
(ARO) grant W911NF-23-1-0026 during this project. D.B. learned about the Carpenter
property from Fedor Sukochev during a visit to UNSW in Sydney, Australia,
in January/February 2025. He would like to thank him and his group
for generously sharing their work on the Carpenter property, and UNSW for
hosting him and providing a very pleasant work environment. The visit to
UNSW was also supported by a grant awarded by the Sydney Mathematical
Research Institute's International Visitor Program.

\section{Preliminaries}\label{sec:prelim}

\subsection{Conventions}
All finite von Neumann algebras carry specified faithful normal tracial
states. For an inclusion $N\subset M$, the trace on $N$ is the restriction
of the trace $\tau$ on $M$, and $\E\colon M\to N$ is the trace-preserving
conditional expectation. We write $\norm{y}$ for the operator norm,
$\norm{y}_1=\tau(\abs y)$ and $\norm{y}_2=\tau(y^*y)^{1/2}$, and
$\norm{y}_{2,A}$ when the trace of an algebra $A$ must be specified;
$\Proj(M)$ is the set of projections of $M$. Inner products are linear
in the first variable; on $L^2(M)$,
$\langle\hat x,\hat y\rangle=\tau(y^*x)$. A corner $eMe$, for
$0\ne e\in\Proj(M)$,
carries the normalized trace $\tau_e(y)=\tau(y)/\tau(e)$.

For self-adjoint $y$, the positive and negative parts $y_+,y_-\ge0$
satisfy $y=y_+-y_-$ and $y_+y_-=0$.

We use repeatedly that in a $\II$ factor $M$, for every nonzero $e\in\Proj(M)$,
there are partial isometries $v_1,\dots,v_m\in M$ with
\begin{equation}\label{eq:partition}
 v_i^*v_i\le e,\qquad \sum_{i=1}^m v_iv_i^*=1 ;
\end{equation}
partition $1$ into finitely many projections of trace at most
$\tau(e)$ and use comparison of projections.

The operator-norm unit ball of $M$ is complete for $\norm{\cdot}_2$
\cite[Prop.~2.6.4]{AP}. Hence a $\norm{\cdot}_2$-Cauchy sequence
$(p_n)$ of projections converges in $\norm{\cdot}_2$ to some $p\in M$
with $\norm p\le1$; then $p=p^*$, since the adjoint is
$\norm{\cdot}_2$-isometric, and $p^2=p$, since
$\norm{p_n^2-p^2}_2\le2\norm{p_n-p}_2$. So the limit is a projection.

\subsection{The basic construction}
Let $N\subset M$ be an inclusion of $\II$ factors. On $L^2(M)$ let
$J\hat x=\widehat{x^*}$, and let $e_N$ be the orthogonal projection onto
$L^2(N)$. The basic construction is $B=\basic MN=(JNJ)'$. It carries a
unique normal semifinite faithful trace $\Tr$ with
$\Tr(xe_Ny)=\tau(xy)$ for $x,y\in M$. For a projection $T\in B$, the
range $TL^2(M)$ is a right $N$-module and
\[
 \Tr(T)=\dim_N\bigl(TL^2(M)\bigr),
\]
the right $N$-dimension; in particular $\Tr(1)=[M:N]$
\cite[Prop.~9.4.2, (9.3), Def.~8.5.4 and Def.~9.4.9]{AP}. Moreover
$M'\cap B=JMJ\cap(JNJ)'=J(N'\cap M)J$.

The index is invariant under compression by projections of the smaller
factor:
\begin{equation}\label{eq:jones-compression}
 [qMq:qNq]=[M:N]\qquad(0\ne q\in\Proj(N)),
\end{equation}
with normalized traces on both corners. For $[M:N]<\infty$ this is due
to Jones \cite{Jones83}. We also use \eqref{eq:jones-compression} when
$[M:N]=\infty$, so we check that $[qMq:qNq]<\infty$ forces
$[M:N]<\infty$. A finite family $m_1,\dots,m_k\in M$ is a
\emph{Pimsner--Popa basis} of $M$ over $N$ if
$x=\sum_jm_j\E(m_j^*x)$ for all $x\in M$, and such a family exists if
and only if $[M:N]<\infty$ \cite{PP86}; see also \cite[Prop.~9.4.8]{AP}.
Since $\E(qyq)=q\E(y)q$, the restriction of $\E$ to $qMq$ is the
trace-preserving expectation onto $qNq$. Suppose $[qMq:qNq]<\infty$,
let $m_1,\dots,m_k$ be a Pimsner--Popa basis of $qMq$ over $qNq$, and
choose $v_1,\dots,v_n\in N$ as in \eqref{eq:partition}, applied to the
$\II$ factor $N$ and $e=q$. For $x\in M$ each $v_i^*xv_l$ lies in $qMq$,
and $v_l\in N$, so
\[
 x=\sum_{i,l}v_i(v_i^*xv_l)v_l^*
 =\sum_{i,j,l}v_im_j\E(m_j^*v_i^*xv_lv_l^*)
 =\sum_{i,j}(v_im_j)\,\E\bigl((v_im_j)^*x\bigr).
\]
Thus the $v_im_j$ form a finite Pimsner--Popa basis of $M$ over $N$,
and $[M:N]<\infty$. Hence \eqref{eq:jones-compression} holds in all
cases, both sides being infinite simultaneously.

\subsection{The Carpenter property}
\begin{definition}[{\cite[Def.~1.1]{Popa89}}]\label{def:Lambda}
For an inclusion $N\subset M$ of $\II$ factors,
\[
 \Lam(M,N)=\{\alpha\in[0,1]:\ \E(p)=\alpha1\text{ for some }p\in\Proj(M)\}.
\]
\end{definition}
If $\E(p)=\alpha1$, trace preservation forces $\tau(p)=\alpha$, and
$\E(1-p)=(1-\alpha)1$. Thus $0,1\in\Lam(M,N)$ and $\Lam(M,N)$ is
invariant under $\alpha\mapsto1-\alpha$.

\begin{lemma}[Pimsner--Popa obstruction]\label{lem:PP}
If $[M:N]=d<\infty$, then $\Lam(M,N)\cap(0,d^{-1})=\emptyset$.
\end{lemma}
\begin{proof}
The Pimsner--Popa inequality gives $\E(a)\ge d^{-1}a$ for every $a\in M_+$
\cite{PP86}; see \cite[1.2.2]{Popa94}. If $\E(p)=\alpha1$ with
$0<\alpha<d^{-1}$, then $p\ne0$ by trace preservation, and compressing
the inequality by $p$ gives $\alpha p\ge d^{-1}p$, a contradiction.
\end{proof}

\begin{definition}\label{def:carpenter}
An inclusion $N\subset M$ of $\II$ factors has the \emph{Carpenter
property} if for every $b\in N$ with $0\le b\le1$ there is a projection
$p\in M$ with $\E(p)=b$.
\end{definition}
The name comes from the Carpenter theorem for masas
\cite{Kadison02a,Kadison02b}. Taking $b=\alpha1$ shows that the
Carpenter property implies $\Lam(M,N)=[0,1]$; by \cref{lem:PP}, it
therefore forces $[M:N]=\infty$.

\subsection{Uniform infinite index}
\begin{definition}\label{def:uniform}
An inclusion $N\subset M$ of $\II$ factors has \emph{uniform infinite
index} if
\[
 [pMp:Np]=\infty\qquad(0\ne p\in\Proj(N'\cap M)).
\]
Here $Np=\{xp:x\in N\}$, with unit $p$ and normalized trace. Since $N$
is a factor and $\E(p)=\tau(p)1$, the map $x\mapsto xp$ is a
trace-preserving isomorphism of $N$ onto $Np$.
\end{definition}
If $N'\cap M=\C1$, uniform infinite index is simply $[M:N]=\infty$.
This condition, for arbitrary subalgebras, is the standing hypothesis
of \cite{Popa19}; Popa notes there the equivalence with condition~(2)
of the next lemma \cite[\S1]{Popa19}, and we include a proof.

\begin{lemma}\label{lem:WM}
For an inclusion $N\subset M$ of $\II$ factors with basic construction
$B=\basic MN$, the following are equivalent:
\begin{enumerate}
\item $N\subset M$ has uniform infinite index;
\item $\Tr(T)=\infty$ for every nonzero projection $T\in M'\cap B$;
\item the bimodule ${}_ML^2(M)_N$ has no nonzero closed subbimodule of
 finite right $N$-dimension.
\end{enumerate}
\end{lemma}
\begin{proof}
(2)$\Leftrightarrow$(3): the closed $M$-$N$ subbimodules of $L^2(M)$ are
exactly the ranges of projections $T\in M'\cap(JNJ)'=M'\cap B$, and the
right $N$-dimension of $TL^2(M)$ is $\Tr(T)$.

(1)$\Leftrightarrow$(2): since $M'\cap B=J(N'\cap M)J$, the nonzero
projections of $M'\cap B$ are the operators $JpJ$ with
$0\ne p\in\Proj(N'\cap M)$, and $JpJ$ is right multiplication by $p$.
Thus $\Tr(JpJ)=\dim_N(L^2(M)p)$, where $L^2(M)p$ is a right $N$-module
because $p$ commutes with $N$. It suffices to show that
$\dim_N(L^2(M)p)<\infty$ if and only if $[pMp:Np]<\infty$.

The subspace $pL^2(M)p\subset L^2(M)p$ is a right $N$-submodule, since
left multiplication by $p$ commutes with the right action of $N$. Up to
rescaling the norm by the constant $\tau(p)^{1/2}$, it is $L^2(pMp)$,
and the right action of $x\in N$ on it is the right action of
$xp\in Np$. Since $x\mapsto xp$ is a trace-preserving isomorphism,
$\dim_N(pL^2(M)p)=\dim_{Np}L^2(pMp)=[pMp:Np]$. Hence
$[pMp:Np]\le\dim_N(L^2(M)p)$. Conversely, choose $v_1,\dots,v_m\in M$ as
in \eqref{eq:partition} with $e=p$. Then $L^2(M)p$ is the orthogonal sum
of the right $N$-submodules $v_iv_i^*L^2(M)p$, and left multiplication
by $v_i$ maps $pL^2(M)p$ onto $v_iv_i^*L^2(M)p$ as a right $N$-module
map. Therefore $\dim_N(L^2(M)p)\le m[pMp:Np]$.
\end{proof}
In the language of Popa's intertwining by bimodules
\cite{Popa06Betti,Popa06Rigidity}, condition~(3) says that no corner of
$M$ can be intertwined into $N$ inside $M$, written $M\nprec_MN$; see
\cite[\S1]{Popa19} and \cite[Def.~6.2]{Jeong26}. Thus uniform infinite
index is exactly $M\nprec_MN$.

\section{A convex correction estimate}\label{sec:correction}

We prove \cref{thm:A} by correcting an approximate lift in two stages. Given a projection $p$ with
$\E(p)\approx b$, this section produces a positive contraction $x$ with
$\E(x)=b$ that
is close to $p$ in an affine sense; \cref{sec:flattening,sec:lifting}
then replace $x$ by a projection.

This section applies to an arbitrary finite tracial inclusion
$N\subset M$. For a projection $p\in M$ define the affine functional
\begin{equation}\label{eq:cost}
 \ell_p(x)=\tau(p)+\tau(x)-2\tau(px),\qquad 0\le x\le1 .
\end{equation}
The trace pairing in this expression is real, and
$\ell_p(x)=\tau(p(1-x))+\tau((1-p)x)\ge0$. More precisely,
\begin{equation}\label{eq:costinterpretation}
 \ell_p(x)=\norm{x-p}_2^2+\tau(x-x^2).
\end{equation}
Thus the cost measures both the distance from $p$ and the failure of $x$
to be a projection. It is affine in $x$, which makes it suitable for
convex duality, and on a projection $q$ it is exactly
$\ell_p(q)=\norm{q-p}_2^2$.

\begin{lemma}[Correction among positive contractions]\label{lem:correction}
Let $p\in M$ be a projection and let $b\in N$ satisfy
$\eta1\le b\le(1-\eta)1$, where $0<\eta\le1/2$. There exists $x\in M$
with
\[
 0\le x\le1,\qquad \E(x)=b,\qquad
 \ell_p(x)\le\eta^{-1/2}\norm{b-\E(p)}_1 .
\]
\end{lemma}
\begin{proof}
Write $E=\E$ and $s=2p-1$. For self-adjoint $y\in M$ we use the
variational formula
\begin{equation}\label{eq:positivepart}
 \tau(y_+)=\max_{0\le z\le1}\tau(yz);
\end{equation}
indeed $\tau(yz)\le\tau(y_+z)\le\tau(y_+)$, and the positive spectral
projection of $y$ gives equality.

\emph{Passing to the dual problem.} The set
$\mathcal K=\{x\in M:0\le x\le1\}$ of positive contractions is weak-*
compact and convex. For
$x\in\mathcal K$ and $h\in N_{\sa}$ define
\[
 L(x,h)=\tau(p)-\tau(sx)+\tau\bigl(h(b-E(x))\bigr).
\]
Regard $N_{\sa}$ as a real normed space. The function $L$ is affine and
continuous in each variable, with the weak-* topology on $\mathcal K$,
so Sion's minimax theorem \cite{Sion58} gives
\[
 \inf_{x\in\mathcal K}\sup_{h\in N_{\sa}}L(x,h)
 =\sup_{h\in N_{\sa}}\inf_{x\in\mathcal K}L(x,h).
\]
If $E(x)\ne b$, the supremum on the left is infinite: take
$h=t(b-E(x))$ and let $t\to+\infty$. If $E(x)=b$, it equals $\ell_p(x)$.
On the right, $\tau(hE(x))=\tau(hx)$, and \eqref{eq:positivepart} applies.
Consequently
\begin{equation}\label{eq:duality}
 \min_{\substack{x\in\mathcal K\\ E(x)=b}}\ell_p(x)
 =\tau(p)+\sup_{h\in N_{\sa}}\bigl\{\tau(hb)-\tau((h+s)_+)\bigr\}.
\end{equation}
The minimum is attained, since its feasible set is weak-* compact and
contains $b$.

\emph{Bounding the dual variables.} If the supremum in
\eqref{eq:duality} can be restricted to an operator-norm ball, then
$\tau(h(b-E(p)))$ is bounded by a constant times $\norm{b-E(p)}_1$. The
spectral gap of $b$ supplies this restriction. Set $\rho=\eta^{-1/2}$ and
\[
 F_b^s(h)=\tau(hb)-\tau((h+s)_+).
\]
We claim that replacing $h$ by $\hat h=\max(-\rho,\min(h,\rho))$ does
not decrease $F_b^s(h)$. Both $h$ and $\hat h$ lie in $N_{\sa}$, and the
argument uses only $\norm s\le1$, so it also applies with $1-b$ and
$-s$.

We first remove the upper spectral tail. Put
\[
 h'=\min(h,\rho),\qquad g=h-h'\ge0,\qquad
 e=\one_{(\rho,\infty)}(h),\qquad A=h'+s .
\]
Then $eh'=\rho e$ and $g=ege$. Let $f=\one_{(-\infty,0]}(A)$, $Q=1-f$
and $C=fAf$. We have $C\le0$ in $fMf$. If $f\ne0$, then $C-\rho f$ is
invertible in that corner with $\norm{(C-\rho f)^{-1}}\le\rho^{-1}$, and
since $Af=fA$,
\[
 ef(C-\rho f)=eAf-\rho ef=esf .
\]
Hence $\norm{ef}\le\rho^{-1}$, and
\begin{equation}\label{eq:angle}
 eQe=e-efe\ge(1-\rho^{-2})e ;
\end{equation}
this also holds if $f=0$. Applying \eqref{eq:positivepart} to $A+g$
with the test projection $Q$ gives
$\tau((A+g)_+)-\tau(A_+)\ge\tau(Qg)$, and by \eqref{eq:angle} and $g=ege$,
\[
 \tau(Qg)=\tau(g^{1/2}Qg^{1/2})\ge(1-\rho^{-2})\tau(g).
\]
On the other hand $\tau(gb)\le(1-\eta)\tau(g)$. Therefore
\[
 F_b^s(h)-F_b^s(h')\le(\rho^{-2}-\eta)\tau(g)=0 .
\]
To remove the lower tail, use $F_b^s(h)=F_{1-b}^{-s}(-h)-\tau(s)$, which
follows from $y_+=y+(-y)_+$, and apply the previous step to $-h'$ with
$1-b$ and $-s$. This replaces $h'$ by $\max(h',-\rho)=\hat h$ without
lowering the objective.

\emph{Estimating the bounded supremum.} By the claim, the supremum in
\eqref{eq:duality} may be taken over $\norm h\le\rho$.
Testing \eqref{eq:positivepart} on the projection $p$, and using $sp=p$,
gives $\tau((h+s)_+)\ge\tau((h+s)p)=\tau(hE(p))+\tau(p)$. Hence
\[
 \min_{\substack{x\in\mathcal K\\ E(x)=b}}\ell_p(x)
 \le\sup_{\substack{h\in N_{\sa}\\ \norm h\le\rho}}\tau\bigl(h(b-E(p))\bigr)
 \le\rho\,\norm{b-E(p)}_1 .
 \qedhere
\]
\end{proof}

\section{Weak mixing and flattening on corners}\label{sec:flattening}

Throughout \cref{sec:flattening,sec:lifting}, $N\subset M$ is an
inclusion of $\II$ factors with uniform infinite index, and
$B=\basic MN$.

If $p$ and $q$ are projections, then
$\norm{q-p}_2^2=\tau(q)+\tau(p)-2\tau(pq)$. Once the trace of $q$ is
prescribed, only one scalar remains to be controlled: the overlap
$\tau(pq)$. This is why we flatten the completely positive (CP) map
$y\mapsto(\E(y),\tau(py))$ rather than $\E$ alone. The construction
takes place inside spectral corners of a positive contraction $x\in M$,
which
generally lie in $M$ rather than in $N$, so we need flattening on every
nonzero corner $eMe$.

\subsection{Flattening of completely positive maps}
\begin{definition}\label{def:weakmixing}
A $P$-$A$ Hilbert bimodule is \emph{left weakly mixing} if it has no
nonzero closed subbimodule of finite right $A$-dimension
\cite{PetersonSinclair12}. A CP map
$\Phi\colon P\to A$ between tracial von Neumann algebras is
\emph{subtracial} if $\tau_A(\Phi(y))\le\tau_P(y)$ for all $y\in P_+$.
\end{definition}

The bimodule $\Hh(\Phi)$ of a normal CP map $\Phi$ is obtained from the
algebraic tensor product $P\odot A$ with the positive semidefinite form
$\langle a\otimes z,c\otimes w\rangle=\tau_A(w^*\Phi(c^*a)z)$, by taking
the quotient by null vectors and completing. The actions are
$x(a\otimes z)y=xa\otimes zy$ for $x\in P$ and $y\in A$.

\begin{theorem}[{\cite[Cor.~1.3]{Jeong26}}]\label{thm:flattening}
Let $P$ be a $\II$ factor and $A$ a tracial von Neumann algebra. If
$\Phi\colon P\to A$ is subtracial and completely positive and $\Hh(\Phi)$
is left weakly mixing, then there is a net $u_i\in\U(P)$ such that
\[
 \norm{\Phi(u_iyu_i^*)-\tau_P(y)\Phi(1_P)}_{2,A}\longrightarrow0
 \qquad(y\in P).
\]
\end{theorem}
No separability is required; if $P$ is separable, the net can be taken
to be a sequence \cite[Cor.~1.3]{Jeong26}. At each application we use
one element $y$ and finitely many tolerances, so a net suffices.

\subsection{Weak mixing on corners}
\begin{lemma}[Extension from a corner]\label{lem:commutant}
Let $P$ be a $\II$ factor contained unitally in a von Neumann algebra
$D$, let $0\ne e\in\Proj(P)$, and choose partial isometries
$v_1,\dots,v_m\in P$ with $v_i^*v_i\le e$ and $\sum_iv_iv_i^*=1$, which is
possible since $P$ is a $\II$ factor (see \eqref{eq:partition}). For
$T\in(ePe)'\cap eDe$ put
$S=\sum_iv_iTv_i^*$. Then $S\in P'\cap D$ and $eSe=T$. If $T$ is a
projection, so is $S$; if moreover $D$ carries a normal semifinite trace
$\Tr$, then $\Tr(S)\le m\Tr(T)$.
\end{lemma}
\begin{proof}
For $a\in P$, the elements $v_j^*av_i$ lie in $ePe$ and commute with
$T$, so
\[
 aS=\sum_{i,j}v_j(v_j^*av_i)Tv_i^*=\sum_{i,j}v_jT(v_j^*av_i)v_i^*
 =\sum_jv_jTv_j^*a=Sa .
\]
Since $ev_i$ and $v_i^*e$ lie in $ePe$,
$eSe=\sum_iev_iTv_i^*e=T\sum_iev_iv_i^*e=T$. If $T$ is a projection,
then, since $v_i^*v_j\in ePe$,
\[
 S^2=\sum_{i,j}v_iT(v_i^*v_j)Tv_j^*=\sum_{i,j}v_iv_i^*v_jTv_j^*=S .
\]
Finally $\Tr(v_iTv_i^*)=\Tr(Tv_i^*v_iT)\le\Tr(T)$.
\end{proof}

\begin{lemma}[Corner weak mixing]\label{lem:uniformcorners}
For every nonzero $e\in\Proj(M)$, the bimodule ${}_{eMe}eL^2(M)_N$ is
left weakly mixing.
\end{lemma}
\begin{proof}
We first identify the subbimodules. The right action of $y\in N$ on
$L^2(M)$ is $\xi\mapsto\xi y=Jy^*J\xi$, so $B=(JNJ)'$ is exactly the
algebra of bounded operators on $L^2(M)$ commuting with the right action
of $N$. Since $e\in M\subset B$, the subspace $eL^2(M)$ is invariant
under this action, and an operator on $eL^2(M)$ commutes with the right
action of $N$ if and only if, extended by $0$ on $(1-e)L^2(M)$, it lies
in $B$, that is, in $eBe$. Consequently, a closed subspace
$K\subset eL^2(M)$ is an $eMe$-$N$ subbimodule if and only if its
orthogonal projection $T$ commutes with the right action of $N$ and with
the left action of $eMe$, that is, $T\in(eMe)'\cap eBe$. Moreover the
right $N$-dimension of $K=TL^2(M)$ is $\Tr(T)$, by \cref{sec:prelim}.

Now suppose such a projection $T\ne0$ had $\Tr(T)<\infty$. \Cref{lem:commutant}, applied to $M\subset B$, gives a
nonzero projection $S\in M'\cap B$ with $\Tr(S)\le m\Tr(T)<\infty$,
contradicting \cref{lem:WM}.
\end{proof}

\begin{lemma}[Compression by a projection of $N$]\label{lem:uniformNcorner}
For every $0\ne q\in\Proj(N)$, the inclusion $qNq\subset qMq$ has uniform
infinite index, and its trace-preserving expectation is the restriction
of $\E$.
\end{lemma}
\begin{proof}
Let $p$ be a nonzero projection in $(qNq)'\cap qMq$. Apply
\cref{lem:commutant} to $N\subset M$ and $e=q$: the element
$s=\sum_iv_ipv_i^*$, with $v_i\in N$, is a nonzero projection in
$N'\cap M$ with $qsq=p$. Since $s$ commutes with $q$, we have $p=qs$, and
the inclusion $(qNq)p\subset pMp$ is the compression of $Ns\subset sMs$
by the projection $qs\in Ns$. By \eqref{eq:jones-compression},
$[pMp:(qNq)p]=[sMs:Ns]=\infty$. Finally $\E(qyq)=q\E(y)q$, and the
restriction of $\E$ to $qMq$ preserves the normalized trace.
\end{proof}

\subsection{Simultaneous flattening}
\begin{proposition}[Simultaneous corner flattening]\label{prop:cornerflat}
Fix a nonzero $e\in\Proj(M)$, a projection $p\in M$ and
$\lambda\in[0,1]$. There is a net of projections $r_i\le e$ with
\begin{equation}\label{eq:flat}
 \tau(r_i)=\lambda\tau(e),\qquad
 \norm{\E(r_i)-\lambda\E(e)}_2\to0,\qquad
 \tau(pr_i)\to\lambda\tau(pe).
\end{equation}
\end{proposition}
\begin{proof}
Equip $P=eMe$ with $\tau_e$, and $A=N\oplus\C$ with
$\tau_A(x,z)=\frac12\tau(x)+\frac12z$. Consider
\[
 \Psi_e\colon P\to A,\qquad \Psi_e(y)=\bigl(\E(y),\tau(py)\bigr).
\]
The second component is a normal positive functional, since
$\tau(py)=\tau(pyp)\ge0$ for $y\ge0$; hence $\Psi_e$ is normal and
completely positive. It is subtracial: for $y\ge0$,
\[
 \tau_A(\Psi_e(y))=\tfrac12\tau(y)+\tfrac12\tau(py)\le\tau(y)
 =\tau(e)\tau_e(y)\le\tau_e(y).
\]

\emph{Splitting the bimodule.} The central projections $(1,0)$ and
$(0,1)$ of $A$ split $\Hh(\Psi_e)$ into two left $P$-invariant summands
with no mixed terms in the inner product. Up to the factor $\frac12$ in
the inner products, which does not affect weak mixing, they are
$\Hh(\E|_P)$ and the GNS space of $\psi(y)=\tau(py)$.

\emph{The $N$-component.} The map $a\otimes x\mapsto ax\in eL^2(M)$,
for $a,c\in P$ and $x,z\in N$, satisfies
\[
 \langle ax,cz\rangle=\tau(z^*c^*ax)=\tau\bigl(z^*\E(c^*a)x\bigr),
\]
so it descends to an isometry from $\Hh(\E|_P)$ into $eL^2(M)$ that
intertwines both actions. Its range is a $P$-$N$ subbimodule of
$eL^2(M)$, so by \cref{lem:uniformcorners} it contains no nonzero
subbimodule of finite right $N$-dimension.

\emph{The scalar component.} For $y\in P=eMe$ we have $y=eye$, so
traciality gives $\psi(y)=\tau(py)=\tau(epe\,y)$. The GNS space of
$\psi$ is the completion of $P$ for the inner product
$\langle a,c\rangle_\psi=\psi(c^*a)$, with $P$ acting by left
multiplication. The map $a\mapsto a(epe)^{1/2}$ from $P$ into
$L^2(P,\tau_e)$ satisfies
\[
 \bigl\langle a(epe)^{1/2},c(epe)^{1/2}\bigr\rangle_{L^2(P,\tau_e)}
 =\tau_e\bigl((epe)^{1/2}c^*a(epe)^{1/2}\bigr)=\tau(e)^{-1}\tau(epe\,c^*a)
 =\tau(e)^{-1}\langle a,c\rangle_\psi
\]
and intertwines left multiplication. Hence, up to the constant factor
$\tau(e)^{-1}$ in the inner product, the GNS space of $\psi$ is the
closure of $P(epe)^{1/2}$ in $L^2(P,\tau_e)$, with $P$ acting by left
multiplication. In particular this representation is normal, being a
subrepresentation of the standard representation of $P$. Right
dimension over $\C$ is Hilbert-space dimension. A nonzero finite-dimensional
invariant subspace would be reducing and would give a nonzero normal
representation of the $\II$ factor $P$ on a finite-dimensional space.
Its kernel would be a weak-* closed ideal, hence zero by factoriality,
which is impossible since $P$ contains arbitrarily many nonzero
orthogonal projections. (If $epe=0$, this summand is zero.)

\emph{Reassembling.} A closed $P$-$A$ subbimodule $K$ splits as
$K(1,0)\oplus K(0,1)$, and
$\dim_AK=\frac12\dim_N(K(1,0))+\frac12\dim_\C(K(0,1))$. So a nonzero
$K$ of finite right dimension would give one in a component, which has
been excluded. Thus $\Hh(\Psi_e)$ is left weakly mixing, and
\cref{thm:flattening} applies.

Choose $r\in\Proj(P)$ with $\tau_e(r)=\lambda$, possible since $P$ is
diffuse, and let $u_i$ be as in \cref{thm:flattening}. Then
$r_i=u_iru_i^*$ satisfies $\tau(r_i)=\lambda\tau(e)$, and the
convergence $\Psi_e(r_i)\to\lambda\Psi_e(e)$ in $L^2(A,\tau_A)$ is
exactly the remaining two statements of \eqref{eq:flat}.
\end{proof}

\section{Approximation by projections and exact lifting}\label{sec:lifting}

We now replace the positive contraction produced by
\cref{lem:correction} by a
projection with nearly the same expectation and nearly the same overlap
with the previous projection. This projection need not be close to the
positive contraction in $\norm{\cdot}_2$; its distance to the previous
projection is
controlled through the affine cost $\ell_p$.

\begin{lemma}[Approximation by projections with overlap control]\label{lem:projection-approx}
For every positive contraction $x\in M$, projection $p\in M$ and
$\varepsilon,\zeta>0$
there is a projection $q\in M$ with
\[
 \tau(q)=\tau(x),\qquad \norm{\E(q)-\E(x)}_2<\varepsilon,\qquad
 \abs{\tau(pq)-\tau(px)}<\zeta .
\]
\end{lemma}
\begin{proof}
Choose $0<\varepsilon_0<\min(\varepsilon/2,\zeta/2)$ and partition $[0,1]$
into finitely many disjoint Borel intervals of diameter at most
$\varepsilon_0$. Let
$e_1,\dots,e_k$ be the nonzero spectral projections of $x$ for these
intervals, and put
\[
 \lambda_j=\frac{\tau(xe_j)}{\tau(e_j)},\qquad x_0=\sum_{j=1}^k\lambda_je_j .
\]
Each $\lambda_j$ lies in the closed convex hull of the corresponding
interval, so $0\le\lambda_j\le1$, $\norm{x-x_0}\le\varepsilon_0$ and
$\tau(x_0)=\tau(x)$. By \cref{prop:cornerflat}, choose projections
$r_j\le e_j$ with $\tau(r_j)=\lambda_j\tau(e_j)$ and
\[
 \norm{\E(r_j)-\lambda_j\E(e_j)}_2<\frac{\varepsilon}{2k},\qquad
 \abs{\tau(pr_j)-\lambda_j\tau(pe_j)}<\frac{\zeta}{2k}.
\]
The orthogonal sum $q=\sum_jr_j$ is a projection of trace
$\sum_j\lambda_j\tau(e_j)=\tau(x_0)=\tau(x)$. Since
$\E(q)-\E(x_0)=\sum_j\bigl(\E(r_j)-\lambda_j\E(e_j)\bigr)$, the triangle
inequality and $L^2$-contractivity of $\E$ give
\begin{align*}
 \norm{\E(q)-\E(x)}_2
 &\le\sum_{j=1}^k\norm{\E(r_j)-\lambda_j\E(e_j)}_2+\norm{\E(x_0-x)}_2\\
 &<k\cdot\frac{\varepsilon}{2k}+\norm{x_0-x}_2
 \le\frac\varepsilon2+\varepsilon_0<\varepsilon .
\end{align*}
Similarly, since $\tau(pq)-\tau(px_0)=\sum_j\bigl(\tau(pr_j)-\lambda_j\tau(pe_j)\bigr)$,
\begin{align*}
 \abs{\tau(pq)-\tau(px)}
 &\le\sum_{j=1}^k\abs{\tau(pr_j)-\lambda_j\tau(pe_j)}+\abs{\tau(p(x_0-x))}\\
 &<k\cdot\frac{\zeta}{2k}+\norm{x_0-x}\le\frac\zeta2+\varepsilon_0<\zeta .
\end{align*}
\end{proof}

\begin{proposition}[Controlled correction]\label{prop:controlled}
Let $b\in N$ with $\eta1\le b\le(1-\eta)1$, $0<\eta\le1/2$. For every
projection $p\in M$ with $\delta=\norm{\E(p)-b}_2>0$ there is a
projection $q\in M$ with
\begin{equation}\label{eq:controlled}
 \norm{\E(q)-b}_2<\frac{\delta}{4},\qquad
 \norm{q-p}_2^2\le(\eta^{-1/2}+1)\delta,\qquad
 \tau(q)=\tau(b).
\end{equation}
\end{proposition}
\begin{proof}
By \cref{lem:correction} there is a positive contraction $x$ with
$\E(x)=b$ and
\[
 \ell_p(x)\le\eta^{-1/2}\norm{b-\E(p)}_1\le\eta^{-1/2}\delta .
\]
Apply
\cref{lem:projection-approx} to $x$ and $p$, with $\varepsilon=\delta/4$
and $\zeta=\delta/2$. The resulting projection $q$ has $\tau(q)=\tau(x)=\tau(b)$
and the required expectation bound, and
\[
 \norm{q-p}_2^2=\tau(q)+\tau(p)-2\tau(pq)=\ell_p(x)+2\tau(p(x-q))
 \le\eta^{-1/2}\delta+2\cdot\frac\delta2 .
 \qedhere
\]
\end{proof}

\begin{corollary}[Targets bounded away from $0$ and $1$]\label{cor:interior}
If $b\in N$ and $\eta1\le b\le(1-\eta)1$ for some $\eta>0$, then
$\E(p)=b$ for some projection $p\in M$.
\end{corollary}
\begin{proof}
Choose any projection $p_0\in M$ and apply \cref{prop:controlled}
inductively; if some $\delta_n=\norm{\E(p_n)-b}_2$ vanishes, stop.
Otherwise
\[
 \delta_n\le4^{-n}\delta_0,\qquad
 \norm{p_{n+1}-p_n}_2\le(\eta^{-1/2}+1)^{1/2}\,2^{-n}\delta_0^{1/2},
\]
so $(p_n)$ is Cauchy in $\norm{\cdot}_2$. By \cref{sec:prelim}, its
limit is a projection $p\in M$, and
$\norm{\E(p)-b}_2\le\norm{p-p_n}_2+\delta_n\to0$.
\end{proof}

\begin{proof}[Proof of \cref{thm:A}]
Let $0\le b\le1$ in $N$, and put $q_1=p_1=\one_{\{1\}}(b)$, so that
$\E(p_1)=q_1=bq_1$. For $k\ge2$ let
$K_k=[1/k,1-1/k]$, and define disjoint Borel sets $I_2=K_2$ and
$I_k=K_k\setminus K_{k-1}$ for $k\ge3$, whose union is $(0,1)$. Let
$q_k=\one_{I_k}(b)$ for $k\ge2$, omitting those that vanish. Then
$\frac1kq_k\le bq_k\le(1-\frac1k)q_k$. By \cref{lem:uniformNcorner}, the
inclusion $q_kNq_k\subset q_kMq_k$ has uniform infinite index and
expectation $\E|_{q_kMq_k}$, so \cref{cor:interior}, applied there,
gives a projection $p_k\le q_k$ with $\E(p_k)=bq_k$. Since the $q_k$ are
mutually orthogonal, the partial sums of $p=\sum_{k\ge1}p_k$ are
increasing projections, so the series converges strongly to a
projection of $M$, and normality of $\E$ gives
\[
 \E(p)=\sum_{k\ge1}bq_k=b .
 \qedhere
\]
\end{proof}

The remaining statements of \cref{thm:A} follow from the remarks after
\cref{def:carpenter,def:uniform}: the Carpenter property implies
$\Lam(M,N)=[0,1]$, and an inclusion with $N'\cap M=\C1$ and
$[M:N]=\infty$ has uniform infinite index.

\section{Unitary bases: crossed products and group algebras}\label{sec:bases}

Uniform infinite index is easy to verify when $L^2(M)$ has an infinite
orthonormal basis over $N$ consisting of unitaries of $M$.

\begin{proposition}\label{prop:unitary-basis}
Let $N\subset M$ be an inclusion of $\II$ factors, and let
$(u_i)_{i\in I}$ be unitaries in $M$ with
\begin{equation}\label{eq:unitary-basis}
 \sum_{i\in I}u_ie_Nu_i^*=1
\end{equation}
strongly on $L^2(M)$. Then $[M:N]=\abs I$. If $I$ is infinite, then
$N\subset M$ has uniform infinite index, and hence the Carpenter
property.
\end{proposition}
\begin{proof}
Each $p_i=u_ie_Nu_i^*$ is a projection in $B=\basic MN$, and by
\eqref{eq:unitary-basis} the $p_i$ are mutually orthogonal. For
$S\in B$, the operator $e_NSe_N$ restricted to $L^2(N)$ commutes with the
right action of $N$, so it is left multiplication by some $a\in N$, and
$e_NSe_N=ae_N$. If $S\ge0$, then $a\ge0$ and
\[
 \Tr(e_NSe_N)=\Tr(ae_N)=\tau(a)=\langle e_NSe_N\hat1,\hat1\rangle
 =\langle S\hat1,\hat1\rangle .
\]
For $T\in B_+$, normality and traciality of $\Tr$ and unitary
invariance under $u_i\in B$ therefore give
\begin{equation}\label{eq:trace-basis}
 \Tr(T)=\sum_{i\in I}\Tr(p_iTp_i)
 =\sum_{i\in I}\Tr\bigl(e_N(u_i^*Tu_i)e_N\bigr)
 =\sum_{i\in I}\langle T\hat u_i,\hat u_i\rangle .
\end{equation}
With $T=1$ this gives $[M:N]=\Tr(1)=\abs I$.

Let $Q\in M'\cap B$ be a nonzero projection and put
$c=\langle Q\hat1,\hat1\rangle=\norm{Q\hat1}^2$. If $c=0$, then
$Q\hat x=xQ\hat1=0$ for all $x\in M$, so $Q=0$; hence $c>0$. Since $Q$
commutes with $u_i\in M$,
$\langle Q\hat u_i,\hat u_i\rangle=\langle u_i^*Qu_i\hat1,\hat1\rangle=c$
for every $i$, and \eqref{eq:trace-basis} gives $\Tr(Q)=c\abs I=\infty$
when $I$ is infinite. By \cref{lem:WM}, $N\subset M$ has uniform
infinite index, and \cref{thm:A} applies.
\end{proof}

Condition \eqref{eq:unitary-basis} says that $(u_i)$ is a Pimsner--Popa
basis of $M$ over $N$; since the $p_i$ are then orthogonal, it is
equivalent to $\E(u_i^*u_j)=\delta_{ij}1$ together with density of the
linear span of $\bigcup_iu_iN$ in $L^2(M)$. Such unitary orthonormal
bases have been studied for inclusions of finite index
\cite{CeccheriniSilberstein06,BakshiGupta22,CrannKribsPereira23}, and
used to produce elements of $\Lam(M,N)$ in \cite{BakshiGuin25};
\cref{prop:unitary-basis} shows that in infinite
index their existence already forces the Carpenter property.

\begin{corollary}[Crossed products]\label{cor:crossed}
Let $N$ be a $\II$ factor and $\sigma$ an action of an infinite discrete
group $G$ on $N$ by trace-preserving automorphisms, and suppose that
$M=N\rtimes_\sigma G$ is a factor. Then $N\subset M$ has the Carpenter
property.
\end{corollary}
\begin{proof}
The canonical unitaries $u_g$, $g\in G$, satisfy
$\E(u_g^*u_h)=\delta_{g,h}1$, and the span of the sets $u_gN$ is dense in
$L^2(M)$. Apply \cref{prop:unitary-basis}.
\end{proof}
If $\sigma$ is outer, then $N'\cap M=\C1$ and \cref{thm:A} applies
directly; \cref{cor:crossed} also covers actions that are not outer,
such as the trivial action of an ICC group, where
$M=N\mathbin{\overline\otimes}L(G)$ and $N'\cap M=L(G)$.

\begin{corollary}[Group algebras]\label{cor:groups}
Let $H\le G$ be ICC groups. Then $L(H)\subset L(G)$ has the Carpenter
property if and only if $[G:H]=\infty$.
\end{corollary}
\begin{proof}
Let $(g_i)_{i\in I}$ be representatives of the left cosets $gH$, and
$u_g$ the canonical unitaries. Then
$E_{L(H)}(u_{g_i}^*u_{g_j})=E_{L(H)}(u_{g_i^{-1}g_j})=\delta_{ij}1$, and
the span of the sets $u_{g_i}L(H)$ contains every $u_g$, so it is dense.
By \cref{prop:unitary-basis}, $[L(G):L(H)]=[G:H]$, and if $[G:H]=\infty$
the inclusion has the Carpenter property. If $[G:H]<\infty$, it fails by
\cref{lem:PP} and the remark after \cref{def:carpenter}.
\end{proof}
\Cref{cor:groups} proves \cite[Conj.~3.15]{HKS}. The relative commutant
$L(H)'\cap L(G)$ plays no role, so the inclusion need not be
irreducible. The case of free groups,
$L(\mathbb F_{S_0})\subset L(\mathbb F_S)$ for a proper subset
$S_0\subset S$ with at least two elements, is \cite[Thm.~3.9]{HKS}.

\section{Examples without uniform infinite index}\label{sec:examples}

The first implication in \eqref{eq:implications} is \cref{thm:A}, and
the second follows from \cref{lem:PP} (see the remark after
\cref{def:carpenter}). For irreducible inclusions the three conditions
are equivalent, since uniform infinite index then reduces to
$[M:N]=\infty$. The two examples below show that neither implication
reverses in general. \Cref{ex:nonuniformyes} has the Carpenter property
without uniform infinite index, and \cref{ex:nonuniformno} has infinite
index without the Carpenter property. Both are built from the
hyperfinite $\II$ factor and contain a corner of index one.

\begin{example}[Carpenter without uniform infinite index]\label{ex:nonuniformyes}
Let $R$ be the hyperfinite $\II$ factor, realized as an infinite tensor
product of copies of $M_2(\C)$. The even tensor positions give a
trace-preserving embedding $\sigma\colon R\to R$, and the odd positions
generate a factor $D\cong R$ with $R=\sigma(R)\mathbin{\overline\otimes}D$.
Set
\[
 M=M_2(R),\qquad
 N=\left\{\begin{pmatrix}x&0\\0&\sigma(x)\end{pmatrix}:x\in R\right\},
\]
with the normalized matrix trace tensored with $\tau_R$ on $M$. Both are
$\II$ factors. With $F=\sigma^{-1}\circ E_{\sigma(R)}$, pairing with
diagonal elements of $N$ shows
\begin{equation}\label{eq:diagexpectation}
 \E\begin{pmatrix}a&v\\w&c\end{pmatrix}
 =\begin{pmatrix}z&0\\0&\sigma(z)\end{pmatrix},\qquad z=\frac{a+F(c)}2 .
\end{equation}
Given a positive contraction $x\in R$, let $q=\one_{[1/2,1]}(x)$ and $y=2x-q$; then $q$
is a projection and $0\le y\le1$. Choose $g\in D$ with uniform spectral
distribution on $[0,1]$ and put $r=\one_{[0,\infty)}(\sigma(y)-g)$; the two
operators in the difference commute. Integrating over the $D$-variable,
\[
 E_{\sigma(R)}(r)=\int_0^1\one_{[t,1]}(\sigma(y))\,dt=\sigma(y),
\]
by the identity $\int_0^1\one_{\{t\le\lambda\}}\,dt=\lambda$
($\lambda\in[0,1]$) applied through spectral
calculus. Thus $F(r)=y$, and by \eqref{eq:diagexpectation} the projection
$p=\operatorname{diag}(q,r)$ satisfies
$\E(p)=\operatorname{diag}(x,\sigma(x))$. Hence $N\subset M$ has the
Carpenter property, and $[M:N]=\infty$ by \cref{lem:PP}. On the other
hand $f=\operatorname{diag}(1,0)$ lies in $N'\cap M$ and $Nf=fMf$, so
$[fMf:Nf]=1$ and the index is not uniformly infinite.
\end{example}

\begin{remark}\label{rem:HKS-example}
\Cref{ex:nonuniformyes} is also covered by \cite[Thm.~3.1]{HKS}, which
gives the Carpenter property for a subfactor whenever $N'\cap M$ has
projections of every trace in $[0,1]$. Here $N'\cap M$ contains $f$ and
the projections $\operatorname{diag}(0,d)$ with $d\in\Proj(D)$, since $D$
commutes with $\sigma(R)$, and their sums have traces
$\frac12(\lambda+\tau_R(d))$, $\lambda\in\{0,1\}$, which fill $[0,1]$.
That theorem is deduced from \cite[Lemma~3.2]{HKS}, which lifts a target
piece by piece along a partition of unity in $N'\cap M$; since
projections of $N'\cap M$ have scalar expectation, no condition on the
corners is needed, and our construction on the corner $1_M-f$ is a
continuous version of the same idea. The relative commutant need not be
diffuse: if $p_n\in\Proj(R)$ have traces $2^{-n}$ and sum $1$, and
$\theta_n\colon R\to p_nRp_n$ are isomorphisms, then
$N=\{\sum_n\theta_n(x):x\in R\}\subset R$ has
$N'\cap R=\overline{\operatorname{span}}\{p_n\}$, so it has the Carpenter
property by \cite[Thm.~3.1]{HKS}, although every corner $Np_n=p_nRp_n$
has index one.
\end{remark}

\begin{example}[Infinite index without Carpenter]\label{ex:nonuniformno}
Realize the hyperfinite $\II$ factor as
$R=\overline\bigotimes_{k\in\Z}M_2(\C)$, with the bilateral Bernoulli
shift $\beta$, and let $Q=R\rtimes_\beta\Z$ with canonical unitary $u$.
The nonzero powers of $\beta$ are mixing, hence outer. If
$x=\sum_nx_nu^n\in R'\cap Q$, then $ax_n=x_n\beta^n(a)$ for all $a\in R$
and $n\in\Z$. For such an $x_n$, both $x_nx_n^*$ and $x_n^*x_n$ commute
with $R$, so they are scalars, and if $x_n\ne0$ then $x_n$ is a multiple
of a unitary implementing $\beta^n$; hence $x_n=0$ for $n\ne0$, while
$x_0\in R'\cap R=\C1$. Thus $R'\cap Q=\C1$. The Fourier decomposition gives $[Q:R]=\infty$. The factor $Q$ is amenable, as the crossed product of
an amenable von Neumann algebra by an amenable group
\cite[Cor.~4.2]{AD79}, and
hence hyperfinite by Connes' theorem \cite{Connes76}; see also
\cite[Thms.~11.1.3 and 11.2.2]{AP}. Fix $e\in\Proj(R)$ with
$\tau_Q(e)=s\in(0,1)$. Since $Q$ and $eRe$ are both hyperfinite, the
uniqueness of the hyperfinite $\II$ factor \cite[Thm.~11.2.2]{AP} gives
a trace-preserving isomorphism $\iota\colon Q\to eRe$ (normalized
traces). Set $1_M=\operatorname{diag}(1,e)$ and
\[
 M=1_MM_2(Q)1_M,\qquad
 N=\left\{\begin{pmatrix}x&0\\0&\iota(x)\end{pmatrix}:x\in Q\right\},
\]
so that $\tau_M(\operatorname{diag}(a,w))=(\tau_Q(a)+\tau_Q(w))/(1+s)$.
With $F=\iota^{-1}\circ E_{eRe}$, the corner expectation taken with
normalized traces,
\begin{equation}\label{eq:weightedexpectation}
 \E\begin{pmatrix}a&v\\w&c\end{pmatrix}
 =\begin{pmatrix}z&0\\0&\iota(z)\end{pmatrix},\qquad
 z=\frac{a+sF(c)}{1+s}.
\end{equation}
Let $f=\operatorname{diag}(1,0)$ and $g=1_M-f$. Then
\[
 \tau_M(f)=t:=\frac1{1+s}>\frac12,\qquad [fMf:Nf]=1,\qquad
 [gMg:Ng]=[eQe:eRe]=\infty,
\]
the last equality by \eqref{eq:jones-compression}. This also forces
$[M:N]=\infty$, since finite right $N$-dimension of $L^2(M)$ would give
finite right $Ng$-dimension to $gL^2(M)g$. The index-one corner shows
that the index is not uniformly infinite. In fact
$N'\cap M=\C f\oplus\C g$: the diagonal entries of a commuting element
are scalar by factoriality of $Q$ and irreducibility of
$eRe\subset eQe$, and a nonzero off-diagonal entry $v\in Qe$ would
satisfy $xv=v\iota(x)$ for all $x\in Q$, forcing $vv^*\in\C1$ and
$v^*v\in\C e$, so that its polar part would have initial projection $e$
and final projection $1$, contradicting $\tau_Q(e)<1$.

We show that $N\subset M$ fails the Carpenter property. Choose
$\gamma\in(1-t,t)$ and suppose $p\in\Proj(M)$ satisfies
$\E(p)=\gamma1_M$. Positivity of the second diagonal entry and
\eqref{eq:weightedexpectation} give $0\le fpf\le\frac\gamma tf$, so
$\norm{fpf}<1$. But $\tau_M(p)=\gamma$, so
\[
 \tau_M(p\wedge f)\ge\tau_M(p)+\tau_M(f)-1=\gamma+t-1>0,
\]
and on the nonzero range of $p\wedge f$ the operator $fpf$ is the
identity, a contradiction. In particular $\Lam(M,N)$ misses the interval
$(1-t,t)$.
\end{example}

In \cref{ex:nonuniformno}, the obstruction is a finite-index corner
carrying a large part of the trace; the mere presence of a finite-index
corner is not enough, as \cref{ex:nonuniformyes} shows.

\section{Further remarks and questions}\label{sec:final}

\begin{remark}[Infinite index]\label{rem:infinite-index}
Under uniform infinite index, $\Lam(M,N)=[0,1]$ by \cref{thm:A}, and
the same holds for $pNp\subset pMp$ for every $0\ne p\in\Proj(N)$, by
\cref{lem:uniformNcorner}. This answers Popa's
\cite[Problems~5.4.1 and~5.4.6]{Popa89} for these inclusions. Without
uniform infinite index, $\Lam(M,N)$ is not determined by the index
(\cref{ex:nonuniformyes,ex:nonuniformno}). It would be interesting to
describe $\Lam(M,N)$ in terms of the corners $Np\subset pMp$,
$p\in\Proj(N'\cap M)$, and their traces.
\end{remark}

\begin{remark}[Extreme points]\label{rem:extreme}
A second route to the Carpenter property goes through extreme points.
For $b\in N$ with $0\le b\le1$, the set
$K_b=\{x\in M:0\le x\le1,\ \E(x)=b\}$ is nonempty, convex and weak-*
compact, so it has extreme points by the Krein--Milman theorem. The
argument of \cite[Lemma~2.1]{DFHS12}, written there for masas, shows
that every extreme point of every $K_b$ is a projection provided
\begin{equation}\label{eq:corner-kernel}
 \ker\E\cap eMe\ne\{0\}\qquad(0\ne e\in\Proj(M)).
\end{equation}
Indeed, if $x\in K_b$ is not a projection, then $e=\one_{[\epsilon,1-\epsilon]}(x)\ne0$
for some $\epsilon>0$, and a self-adjoint $h\in\ker\E\cap eMe$ with
$\norm h\le\epsilon$ gives $x\pm h\in K_b$. Thus \eqref{eq:corner-kernel}
implies the Carpenter property; this is the noncommutative
Lyapunov--Knowles theorem \cite{AzarniaWright82}. Condition
\eqref{eq:corner-kernel} is verified in \cite[Prop.~3.10]{HKS} when $M$
has nonseparable predual and $N$ has separable predual, and it holds
for every $e$ that dominates a
nonzero projection $q\in N$, provided $N\ne M$: for $x\in qMq$, the
element $x-\E(x)$ lies in $qMq$, and $qMq\ne qNq$ by
\eqref{eq:jones-compression}.

Condition \eqref{eq:corner-kernel} is not necessary. If $\E$ is
injective on $eMe$, then $e/2$ is an extreme point of $K_{\E(e)/2}$,
since $e/2=(a+c)/2$ with $a,c\in K_{\E(e)/2}$ forces $a,c\in eMe$ and
$\E(a-e/2)=0$; and in \cref{ex:nonuniformyes}, which has the Carpenter
property, $\E$ is injective on $fMf$. We do not know whether uniform
infinite index implies \eqref{eq:corner-kernel}. The proof of
\cref{thm:A} avoids this question: flattening only shows that $\E$ is
not bounded below on $eMe$, and the exact projection comes from the
controlled iteration of \cref{sec:lifting}.
\end{remark}

\begin{remark}[Orthogonal subalgebras]\label{rem:orthogonal}
If a von Neumann subalgebra $K\subset M$ is orthogonal to $N$, that is,
$\E(K)\subseteq\C1$ \cite{Popa83}, then every projection $p\in K$ satisfies
$\E(p)=\tau(p)1$. It is conjectured in \cite[Conj.~3.8]{HKS} that a
subfactor $N\subset M$ has the Carpenter property if and only if some
$K\subset M$ orthogonal to $N$ has projections of every trace in
$[0,1]$. By \cref{thm:A}, one direction of this
conjecture would give, for every inclusion of uniform infinite index, a
diffuse subalgebra of $M$ orthogonal to $N$. After passing to the
ultrapower such subalgebras exist: every separable subalgebra of
$M^\omega$ is unitarily conjugate to one orthogonal to $N^\omega$
\cite[Cor.~1.2]{Popa19}.
Constructing one in $M$ requires commuting families of projections whose
products all have scalar expectation, which \cref{thm:A} does not
provide.
\end{remark}

\bibliographystyle{amsplain}
\bibliography{refs}

\end{document}